\documentclass[11pt,reqno]{amsart}
\usepackage{amssymb}
\usepackage{amsfonts}
\usepackage{amsmath}
\usepackage{graphicx}
\usepackage{hyperref}
\usepackage{float}
\usepackage{epstopdf}
\usepackage[usestackEOL]{stackengine}
\usepackage{color}
\usepackage{comment}
\usepackage{bm}

\allowdisplaybreaks
\newtheorem{theorem}{Theorem}[section]
\newtheorem{proposition}[theorem]{Proposition}
\newtheorem{lemma}[theorem]{Lemma}
\newtheorem{corollary}[theorem]{Corollary}

\newtheorem{definition}[theorem]{Definition}

\newenvironment{notation}{\smallskip{\sc Notation.}\rm}{\smallskip}

\def\SG{\mathcal{SG}}

\numberwithin{equation}{section}

\begin{document}
\title[A trace theorem for Sobolev spaces on the Sierpinski gasket]{A trace theorem for Sobolev spaces on the Sierpinski gasket}

\author{Shiping Cao}
\address{Department of Mathematics, Cornell University, Ithaca 14853, USA}
\email{sc2873@cornell.edu}
\thanks{}

\author{Shuangping Li}
\address{\Longstack[l]{Program in Applied and Computational Mathematics, Princeton University, NJ 08544-1000, USA}}
\email{sl31@princeton.edu}
\thanks{}

\author{Robert S. Strichartz}
\address{Department of Mathematics, Cornell University, Ithaca 14853, USA}
\email{str@cornell.math.edu}

\author{Prem Talwai}
\address{Department of Mathematics, Cornell University, Ithaca 14853, USA}
\email{pmt55@cornell.edu}

\date{}


\begin{abstract}
We give a discrete characterization of the trace of a class of Sobolev spaces on the Sierpinski gasket to the bottom line. This includes the $L^2$ domain of the Laplacian as a special case. In addition, for Sobolev spaces of low orders, including the domain of the Dirichlet form, the trace spaces are Besov spaces on the line. 
\end{abstract}
\maketitle

\section{Introduction}\label{intro}
This work deals with the restriction problem for functions in Sobolev spaces on the Sierpinski gasket ($\SG$) to the bottom line.  A special case was studied by A. Jonsson in \cite{jonsson1}, where the trace for the Dirichlet form was characterized. 

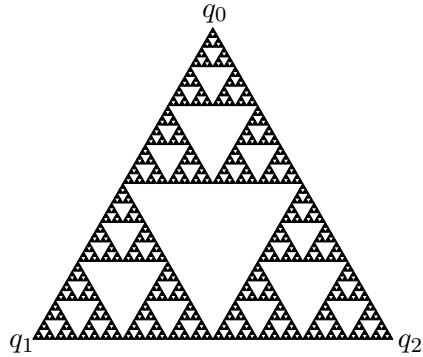
\begin{figure}[h]
	\includegraphics[width=5cm]{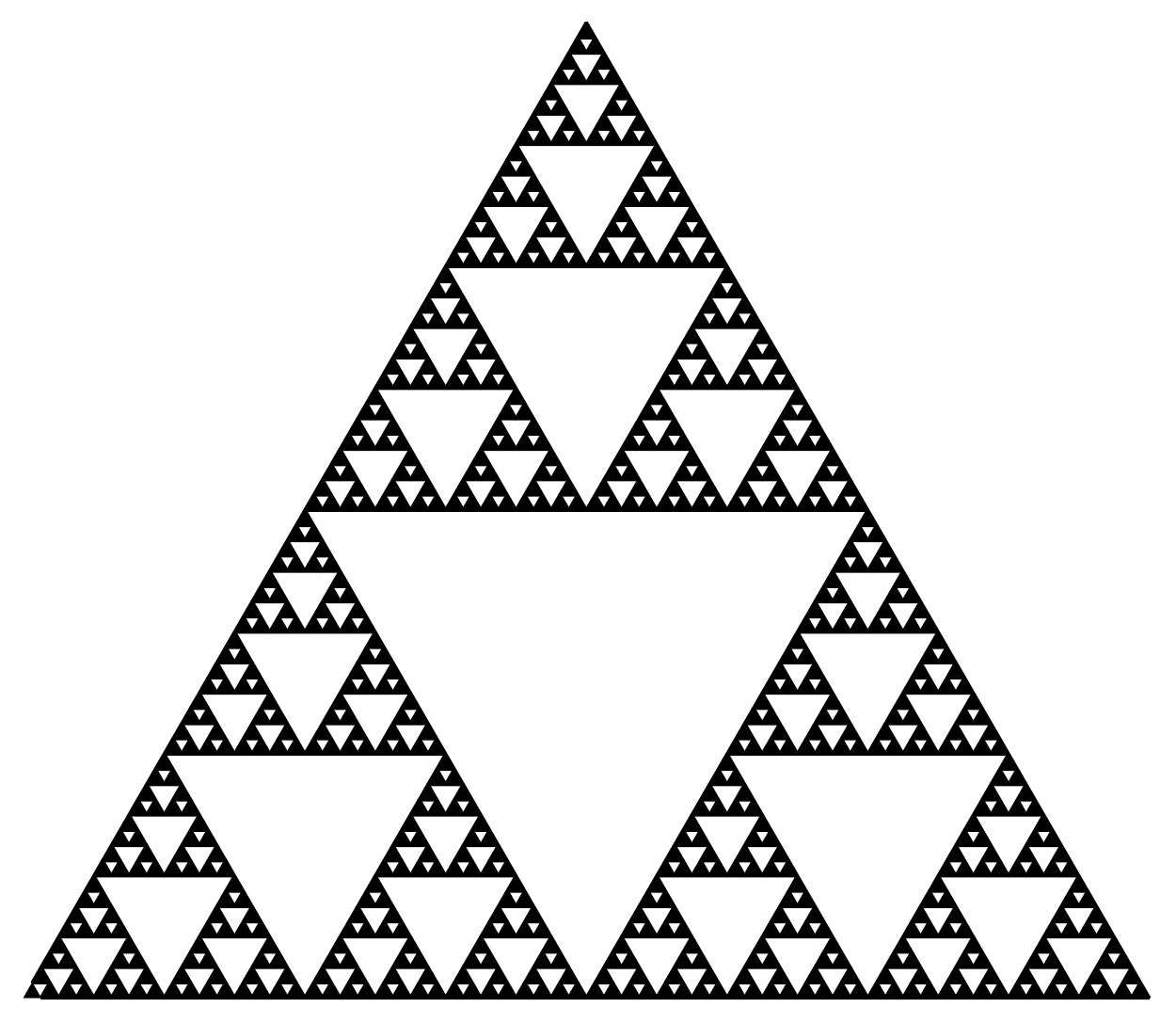}
	\begin{picture}(0,0)
    \put(-5,0){$q_2$}
    \put(-152,0){$q_1$}
    \put(-79,125){$q_0$}
	\end{picture}	
	\caption{the Sierpinski gasket.}\label{SG}
\end{figure}

Let's briefly review Jonsson's result here. $\SG$ is the attractor of the iterated function system
(i.f.s.) in the plane $\mathbb{R}^2$
$$F_i(x)=\frac{1}{2}x+\frac{q_i}{2},\quad i=0,1,2,$$
where $q_0,q_1,q_2$ are vertices of an equilateral triangle. See Figure \ref{SG} for a picture of $\SG$. We identify $I=[0,1]$ with the bottom line $\overline{q_1q_2}$ by  
\begin{equation}\label{eqn11}
t=tq_1+(1-t)q_2,\quad \forall t\in I.
\end{equation}
Consider the standard Dirichlet form $(\mathcal{E},dom\mathcal{E})$ on $\SG$. This form was studied in connection with the Brownian motion on $\SG$ (see \cite{BP,jonsson2,Lindstrom}), and was constructed in a pure analytically approach by J. Kigami in \cite{k1,k2}. A. Jonsson showed the following theorem.

\begin{theorem}[A. Jonsson]\label{thm11}
	Let $\alpha=\frac{\log2+\log5-\log3}{2\log2}(\approx0.868483)$ and $B_\alpha^{2,2}(I)$ be the Besov space on $I$. Then $dom\mathcal{E}|_{I}=B_\alpha^{2,2}(I)$. 
\end{theorem}

The above theorem was extended to a wide class of self-similar sets \cite{HK}, where the trace theorem for Dirichlet forms to self-similar subsets were established. Also, read \cite{HK,kumagai} for an application of the trace theorems to penetrating processes.

Recently, many related works emerge, including the trace theorem on the middle line of $\SG$ (see \cite{cq4,bvphalf}), and boundary value problems on the upper half domain of $\SG$ (see \cite{cq2,bvp2,Hua,bvp1}). However, there have not been further results telling us what is the trace of the domain of the Lapalcian and other Sobolev spaces on $\SG$ to $I$. In this work, we will give an answer to the above question. 

Below, we briefly introduce our results. We choose $\mu$ to be the Hausdoff measure on $\SG$, satisfying $\mu(\SG)=1$ and $\mu(A)=\sum_{i=0}^2\frac{1}{3}\mu(F_i^{-1}A)$ for each Borel set $A$. For $u\in dom\mathcal{E}$, we say $u\in dom_{L^2}\Delta(\SG)$ with $\Delta u=f$ if 
\[\mathcal{E}(u,v)=-\int_{\SG} fvd\mu\]
holds for each $v\in dom_0\mathcal{E}:=\{v\in dom\mathcal{E}:v(q_i)=0,\text{ for } i=0,1,2\}$. In our work, we consider Sobolev spaces $L^2_\sigma(\SG)$ with $0\leq\sigma\leq 2$, which can be defined as follows.

\begin{definition}\label{def12}
	Define the Sobolev space $L^2_0(\SG)=L^2(\SG)$ with norm $\|u\|_{L^2_0(\SG)}=\|u\|_{L^2(\SG)}$, and define $L^2_2(\SG)=dom_{L^2}\Delta(\SG)$ with norm $\|u\|^2_{L^2_2(\SG)}=\|u\|^2_{L^2(SG)}+\|\Delta u\|^2_{L^2(SG)}$.
	
	For $0\leq \sigma\leq 2$, we define Sobolev spaces to be $L^2_\sigma(\SG)=[L^2_0(\SG),L^2_2(\SG)]_{\sigma/2}$, where $[X,Y]_\theta$ denotes the complex interpolation space of $X,Y$.
\end{definition}

In \cite{s2}, Strichartz gave systematic discussions on Sobolev spaces and other function spaces, where Sobolev spaces were defined in more general settings. Related works on properites of Sobolev spaces can be found in \cite{cq1,cq4,GL,HW,HM,pseudo}. Also read \cite{HKM,sw} for recent developments on $p-$Laplacian and the corresponding $L^p$ Sobolev spaces.

In our first main result, we have two critical orders 
\[b_1=\frac{\log3}{\log5}(\approx 0.682606),\quad  b_2=\frac{\log\frac{25(17-\sqrt{73})}{36}}{\log5}(\approx 1.09991).\]
Also, we define the function $\alpha(\sigma):=\frac{\log2+\sigma\log5-\log3}{2\log2}(\approx1.16096\sigma-0.292481)$, which is the unique number such that $5^{\sigma}3^{-1}=2^{2\alpha-1}$. Noticing that $L^2_1(\SG)=dom\mathcal{E}$ (see \cite{cq4}) is included as a special case, the following Theorem \ref{thm13} can be viewed as a direct extension of Jonsson's theorem.

 \begin{theorem}\label{thm13}
 	Let $b_1<\sigma<b_2$ and $\alpha=\alpha(\sigma)$. Then $L^2_\sigma(\SG)|_{I}=B^{2,2}_\alpha(I)$.
 \end{theorem}
 
 $b_1$ in the above theorem is the critical order for the continuity of functions in Sobolev spaces (see \cite{cq1,GL,s2}), and one can check $\alpha(b_1)=\frac{1}{2}$, the well known critical order for Sobolev spaces on the line. The complicated upper bound $b_2$ has an explanation in Corollary \ref{coro33}, where the trace of harmonic functions is studied. 
 
On the other hand, Besov spaces on the line segment are no longer the trace spaces of Sobolev spaces for higher orders. To describe the trace spaces, we will define a difference operator $D$. To be more precise, we define $Df(n,k)$ to be a linear combination of the values of $f\in C(I)$ at $\frac{k}{2^n}$ and some neighbouring points. The space $\mathcal{T}_\sigma$ will be discretely characterized as
\[
\mathcal{T}_{\sigma}=\big\{f\in C(I):\sum_{n=2}^\infty\sum_{k=1}^{2^n-1} 5^{\sigma n}3^{-n}|Df(n,k)|^2<\infty\big\}.
\]
Details can be found in Definition \ref{def35} and \ref{def45}. We will prove the following trace theorem.

\begin{theorem}\label{thm14}
	Let $b_1<\sigma\leq2$. Then $L^2_\sigma(\SG)|_{I}=\mathcal{T}_\sigma$.
\end{theorem}

As supplement, we will show that $B^{2,2}_\alpha(I)\subset \mathcal{T}_\sigma$ for  $b_1<\sigma<\frac{\log6}{\log5}(\approx 1.11328)$, and $B^{2,2}_\alpha(I)= \mathcal{T}_\sigma$ for $b_1<\sigma<b_2$. In addition, $\mathcal{T}_\sigma$ is stable under complex interpolation.

In the end, we briefly introduce the structure of this paper. In section 2, we will review the Dirichlet form and harmonic functions on $\SG$, and introduce some notations and tools. In section 3, we will prove Theorem \ref{thm13}. Some preparations for Theorem \ref{thm14} will be included. In section 4, we will construct the trace space $\mathcal{T}_\sigma$, and prove Theorem \ref{thm14}. In Section 5, we will talk about some related results.

Throughout the paper, we always use the notation $f\lesssim g$ if there is a constant $C>0$ such that $f\leq Cg$, and write $f\asymp g$ if $f\lesssim g$ and $g\lesssim f$. Also, we will keep using the critical numbers $b_1,b_2$ and the function $\alpha(\sigma)$ without further specifying.

\section{The Dirichlet form and harmonic functions}
For convenience of readers, we briefly reivew the Dirichlet form and the harmonic functions on $\SG$ in this section. Some easy lemmas and important tools will also be given. More details can be found in books \cite{k3,s1}.

Recall that $\SG$ is the attractor of the i.f.s $\{F_i\}_{i=0}^2$, i.e.
$$\SG=\bigcup_{i=0}^2F_i\SG.$$
We call each $F_i\SG$ a level-$1$ cell. More generally, define $W_n=\{0,1,2\}^n$ for $n\geq 1$, and set $W_0=\{\emptyset\}$ for uniformity. For each finite word $w=w_1w_2\cdots w_n\in W_*=\bigcup_{m=0}^\infty W_m$, we denote $|w|=n$ the length of the word, and write $F_w=F_{w_1}F_{w_2}\cdots F_{w_n}$ for short. In particular, $F_\emptyset=Id$ is the identity map. We call $F_w\SG$ a level-$n$ cell if $|w|=n$.

We call $V_0=\{q_0,q_1,q_2\}$ the set of boundary vertices of $\SG$, and define the set of level-n vertices $V_n=\bigcup_{w\in W_n} F_wV_0$. For convenience, for $n\geq 1$, let 
\[\tilde{V}_n=V_n\setminus V_{n-1},\]
and set $\tilde{V}_0=V_0$. The set of vertices $V_*=\bigcup_{n=0}^\infty V_m=\bigsqcup_{n=0}^\infty \tilde{V}_n$ is a dense subset of $\SG$.

On $\SG$, J. Kigami \cite{k1,k2} constructed the self-similar energy form by defining it as the limit of a sequence of discrete Dirichlet forms on $V_n$. For each $u\in l(V_n)$, define 
 \[
 \mathcal{E}_n(u)=(\frac{5}{3})^n\sum_{w\in W_n}\sum_{i\neq j}(u(F_wq_i)-u(F_wq_j))^2.
 \]
$\{\mathcal{E}_n(u)\}$ is a nondecreasing sequence for each $u\in C(\SG)$, so we can define $\mathcal{E}(u)=\lim\limits_{n\to\infty} \mathcal{E}_n(u)$. Set $dom\mathcal{E}=\{u\in C(\SG):\mathcal{E}(u)<\infty\}$. For $u,v\in dom\mathcal{E}$, we can use polarization to give a bilinear form
 \[\mathcal{E}(u,v)=\frac{1}{4}\big(\mathcal{E}(u+v)-\mathcal{E}(u-v)\big).\]
It is well known that $(\mathcal{E},dom\mathcal{E})$ is a local regular Dirichlet form on $\SG$ with the Hausdoff measure $\mu$.


	


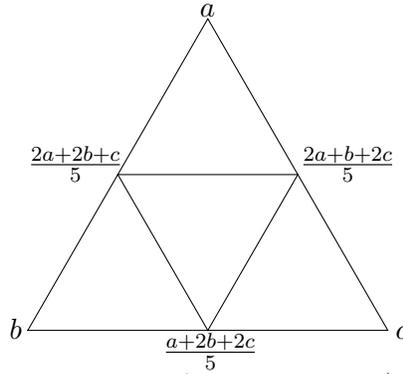
\begin{figure}[h]
	\includegraphics[width=5cm]{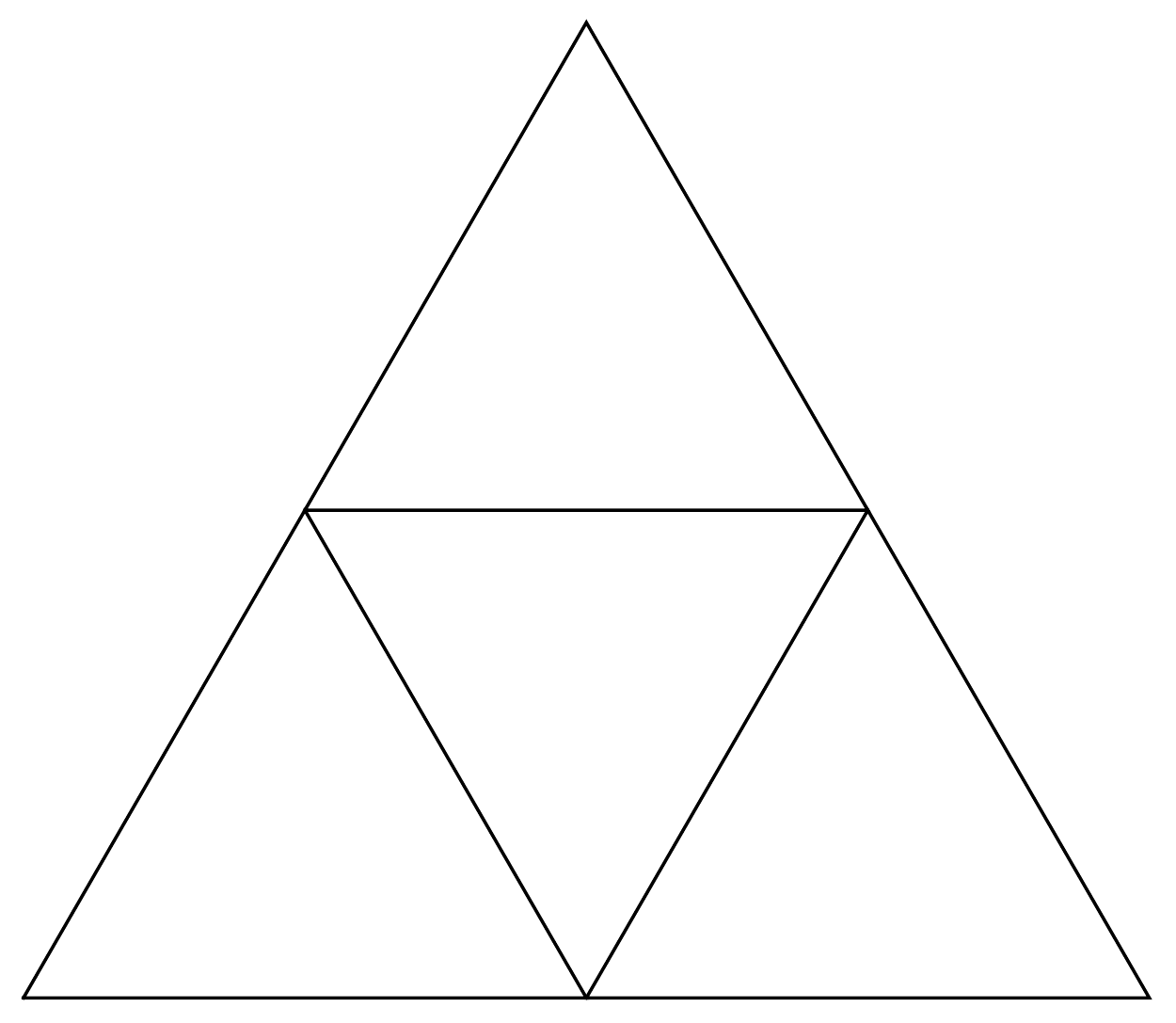}
	\begin{picture}(0,0)
	\put(-78,122){$a$}
	\put(-150,0){$b$}
	\put(-4,0){$c$}
	\put(-143,63){$\frac{2a+2b+c}{5}$}
    \put(-40,63){$\frac{2a+b+2c}{5}$}
     \put(-92,-8){$\frac{a+2b+2c}{5}$}
	\end{picture}
	\caption{The harmonic function with $h(q_0)=a,h(q_1)=b,h(q_2)=c$.}\label{extensionalgorithm}
\end{figure}

Given any boundary value $h_0\in l(V_0)$, there is a unique extension $h_1\in l(V_1)$ that minimizes the energy $\mathcal{E}_1$, i.e. $\mathcal{E}_1(h_1)=\min\{\mathcal{E}_1(u):u\in l(V_1),u|_{V_0}=h_0\}$. The extension algorithm is shown in Figure \ref{extensionalgorithm}.

The above algorithm is local, which means it can be applied to each cell. So we get a sequence of extensions $h_n\in l(V_n)$ that minimize the energy $\mathcal{E}_n$, and $h_n$ converges to $h\in dom\mathcal{E}$. See \cite{s1} for details. $h$ is called a harmonic function, and we denote by $\mathcal{H}_0$ the space of harmonic functions on $\SG$. Clearly, $\mathcal{H}_0$ is of three dimension, since each harmonic function is uniquely determined by its boundary values. The following lemma can be derived from direct computation.

\begin{lemma}\label{lemma21}Let $h$ be a harmonic funciton on $\SG$. For $n\geq 2$ and $1\leq k\leq 2^{n-1}$, we have
	\[h(\frac{2k-1}{2^n})=\begin{cases}
	\frac{4}{5}h(\frac{k}{2^n})+\frac{8}{25}h(\frac{k-1}{2^n})-\frac{3}{25}h(\frac{k+1}{2^n}),\text{ if $k$ is odd},\\
	\frac{4}{5}h(\frac{k-1}{2^n})+\frac{8}{25}h(\frac{k}{2^n})-\frac{3}{25}h(\frac{k-2}{2^n}),\text{ if $k$ is even}.
	\end{cases}\]
\end{lemma}
\textit{Proof.} By using the harmonic extension algorithm twice, we have
\[\begin{aligned}
h(\frac{1}{2})=\frac{2}{5}h(0)+&\frac{2}{5}h(1)+\frac{1}{5}h(q_0),\\
h(\frac{1}{4})=\frac{16}{25}h(0)+\frac{1}{5}h(1)+\frac{4}{25}h(q_0),&\quad h(\frac{3}{4})=\frac{1}{5}h(0)+\frac{16}{25}h(1)+\frac{4}{25}h(q_0),
\end{aligned}\]
where $0=q_1,1=q_2,\frac{1}{2}=F_1q_2,\frac{1}{4}=F^2_1q_2$ and $\frac{3}{4}=F_2^2q_1$ as we set in equation (\ref{eqn11}). For larger $n$, we can do the same computation locally on the level-$(n-2)$ cell containing $\frac{2k-1}{2^n}$. Then, it is direct to check the lemma. \hfill$\square$\vspace{0.2cm}

Analogously to the definition of harmonic functions, for each $m\geq 1$ and $x\in \tilde{V}_m=V_m\setminus V_{m-1}$, we can define the tent function $\varphi_x$ by giving the initial value on $V_m$ as follows  
\[\varphi_x(y)=\begin{cases}
0,\text{ if }y\neq x,\\
1,\text{ if }y=x,
\end{cases}\]
and taking harmonic extension in $\SG\setminus V_m$. Clearly, $\varphi_x$ is harmonic in each level-$m$ cell, and Lemma \ref{lemma21} holds for $\varphi_x$ when $n\geq m+2$.

In our work, we will use the following characterization of Sobolev spaces. For the full version and proof, see Theorem 7.11 in \cite{cq1}.

\begin{theorem}\label{th22}
	Let $b_1=\frac{\log3}{\log5}<\sigma<2-\frac{\log3}{\log5}$, the series
	$f=h+\sum_{n=1}^\infty\sum_{x\in \tilde{V}_n} c_x\varphi_x$
	belongs to $L^2_\sigma(SG)$ if and only if
	\[\sum_{n=1}^\infty \sum_{x\in\tilde{V}_n} 5^{n\sigma}3^{-n}|c_x|^2<\infty.\]
	In addition, each $f\in L^2_\sigma(SG)$ has a unqiue expanison of the above form, with $\|f\|_{L^2_\sigma(SG)}\asymp (\|h\|^2_{L^2(SG)}+\sum_{n=1}^\infty\sum_{x\in \tilde{V}_n} 5^{n\sigma}3^{-n}|c_x|^2)^{\frac{1}{2}}$.  
\end{theorem}

\section{An extension of A. Jonsson's Theorem}




In this section, we study the trace theorem for Sobolev spaces of low orders. The result, Theorem \ref{thm13}, is a direct extension of A. Jonsson's trace theorem. In the following, we will study the restriction map and the extension map seperately. The two parts together imply Theorem \ref{thm13}.

\subsection{A restriction theorem} In this part, we follow A. Jonsson's idea to show a restriction theorem. Recall the fact from \cite{Kamont} that for $\frac{1}{2}<\alpha<1$, a function $f$ belong to $B^{2,2}_\alpha(I)$ if and only if the following expression 
\[
(|f(0)|^2+|f(1)|^2)^{1/2}+\big(\sum_{n=1}^\infty 2^{2n\alpha}2^{-n}\sum_{k=1}^{2^n}|f(\frac{k}{2^n})-f(\frac{k-1}{2^n})|^2\big)^{1/2}
\]
is finite and the norm of $f$ in $B^{2,2}_\alpha(I)$ is equivalent to this expression.

We introduce the following notation to shorten the above expression. 
\begin{definition}\label{def31}
	Let $f\in C(I)$. Define $A_n(f)$ to be a vector of length $2^n$, such that
	\[A_n(f)_k=f(\frac{k}{2^n})-f(\frac{k-1}{2^n}),\forall 1\leq k\leq 2^n.\] 	
\end{definition}

With the above notation, we have 
\[
\|f\|_{B^{2,2}_\alpha(I)}\asymp \big(|f(0)|^2+|f(1)|^2+\sum_{n=1}^\infty 2^{2n\alpha}2^{-n}\|A_n(f)\|^2_{l^2}\big)^{1/2}.
\]

We begin with harmonic functions. 
\begin{proposition}\label{prop32}
Let $h$ be a harmonic function on $\SG$. Then for $n\geq 0$, we have
\begin{equation}\label{eqn31}
	\|A_{n+2}(h|_I)\|^2_{l^2}=\frac{17}{25}\|A_{n+1}(h|_I)\|^2_{l^2}-\frac{54}{625}\|A_{n}(h|_I)\|^2_{l^2}.
\end{equation}
As a consequence, there exist constants $C_1,C_2$ such that
	\[\|A_n(h|_I)\|^2_{l^2}=C_1(\frac{17+\sqrt{73}}{50})^n+C_2(\frac{17-\sqrt{73}}{50})^n.\]
\end{proposition}
\textit{Proof.} By direct computation and using Lemma \ref{lemma21}, we can verify
\[\begin{aligned}
\|A_2(h|I)\|^2_{l^2}&=\sum_{k=1}^{4}\big(h(\frac{k}{4})-h(\frac{k-1}{4})\big)^2\\
&=\frac{17}{25}\Big(\big(h(\frac{1}{2})-h(0)\big)^2+\big(h(1)-h(\frac{1}{2})\big)^2\Big)-\frac{54}{625}\big(h(1)-h(0)\big)^2\\
&=\frac{17}{25}\|A_{1}(h|_I)\|^2_{l^2}-\frac{54}{625}\|A_{0}(h|_I)\|^2_{l^2}.
\end{aligned}\]
This shows (\ref{eqn31}) for $n=0$. For larger $n$, we can do the same computation locally on each $n$ cell and add up to get (\ref{eqn31}).

The second half of the proposition directly follows (\ref{eqn31}), where $\frac{17\pm\sqrt{73}}{50}$ are zeros of the polynomial $x^2-\frac{17}{25}x+\frac{54}{625}$.\hfill$\square$\vspace{0.2cm}

The critical order $b_2$ introduced before Theorem \ref{thm13} is the solution of the equation $5^{b_2}3^{-1}=(\frac{17+\sqrt{73}}{50})^{-1}$. Noticing that $2^{2\alpha(b_2)-1}=5^{b_2}3^{-1}$, we have the following Corollary.

\begin{corollary}\label{coro33}
	Let $h$ be a harmonic funciton on $\SG$. Then $h|_I\in B^{2,2}_\alpha(I)$ if and only if $\alpha<\alpha(b_2)=\frac{\log\frac{25(17-\sqrt{73})}{54}}{2\log2}(\approx 0.984472)$.
\end{corollary}

Using Proposition \ref{prop32} and Theorem \ref{th22}, we can prove the following restriction theorem. 
\begin{theorem}\label{thm34}
Let $b_1<\sigma<b_2$ and $\alpha=\alpha(\sigma)$. Then, the restriction map $u$ to $u|_I$ is continuous $L^2_\sigma(\SG)\to B^{2,2}_\alpha(I)$.
\end{theorem}

\textit{Proof.} By Theorem \ref{th22}, each $u\in L^2_\sigma(\SG)$ admits a unique expansion $$u=h+\sum_{n=1}^\infty\sum_{x\in \tilde{V}_n}  c_x\varphi_x=\sum_{n=0}^{\infty} \psi_n,$$ 
where we write $\psi_n=\sum_{x\in \tilde{V}_n}\varphi_x, n\geq 1$ and $\psi_0=h$ for convenience. Write $C_n=(\sum_{x\in \tilde{V}_n}c^2_x)^{1/2}$ for short.  

Then, obviously  $A_m(\psi_n|_I)=0$ if $m<n$ and $\|A_n(\psi_n|_I)\|_{l^2}\lesssim C_n$. In addition, by Proposition \ref{prop32}, for $m>n$
\[\|A_m(\psi_n|_I)\|_{l^2}\lesssim \lambda^{m-n}\|A_{n+1}(\psi_n|_I)\|_{l^2}\lesssim \lambda^{m-n}C_n,\]
where $\lambda=5^{-b_2/2}3^{1/2}=2^{-\alpha(b_2)+1/2}$. Then, we have the estimate
\[\begin{aligned}
&(\sum_{m=1}^\infty 2^{2m(\alpha-1/2)}\|A_m(u|_I)\|^2_{l^2})^{1/2}
=\big\|2^{m(\alpha-1/2)}\|A_m(u|_I)\|_{l^2}\big\|_{l^2}\\
\leq& \big\|2^{m(\alpha-1/2)}\sum_{n=0}^m\|A_m(\psi_n|_I)\|_{l^2}\big\|_{l^2}\lesssim \|2^{m(\alpha-1/2)}\sum_{n=0}^m \lambda^{m-n}C_n\|_{l^2}\\
=&\|2^{m(\alpha-1/2)}\sum_{n=0}^m \lambda^{n}C_{m-n}\|_{l^2}\leq \sum_{n=0}^\infty (2^{\alpha-1/2}\lambda)^n\|2^{m(\alpha-1/2)}C_m\|_{l^2}\\
\leq& \big(\sum_{n=0}^\infty (2^{\alpha-1/2}\lambda)^n\big)\big(\|h\|^2_{L^2(SG)}+\sum_{n=1}^\infty\sum_{x\in \tilde{V}_n} 5^{n\sigma}3^{-n}|c_x|^2\big)^{\frac{1}{2}}\lesssim \|u\|_{L^2_\sigma(\SG)},
\end{aligned}\]
where we use Theorem \ref{th22} in the last step. The theorem then follows. \hfill$\square$

\subsection{An extension theorem} In the rest of this section, we develop an extension map as the right inverse of the restriction map. It suffices to modify A. Jonsson's idea. However, we provide another extension map here, as preparation for further developments in Section 4. 

We introduce some new notations here. 

\begin{definition}\label{def35}
	Let $f\in C(I)$. For $n\geq 1$ and $1\leq k\leq 2^n$, define $\tilde{D}f(n,k)$  as 
	\[\tilde{D}f(n,k)=\begin{cases}
	f(\frac{2k-1}{2^{n+1}})-\frac{4}{5}f(\frac{k}{2^n})-\frac{8}{25}f(\frac{k-1}{2^n})+\frac{3}{25}f(\frac{k+1}{2^n}),\text{ if $k$ is odd},\\
	f(\frac{2k-1}{2^{n+1}})-\frac{4}{5}f(\frac{k-1}{2^n})-\frac{8}{25}f(\frac{k}{2^n})+\frac{3}{25}f(\frac{k-2}{2^n}),\text{ if $k$ is even}.
	\end{cases}\]
\end{definition}

\begin{figure}[h]
	\includegraphics[width=5.92cm]{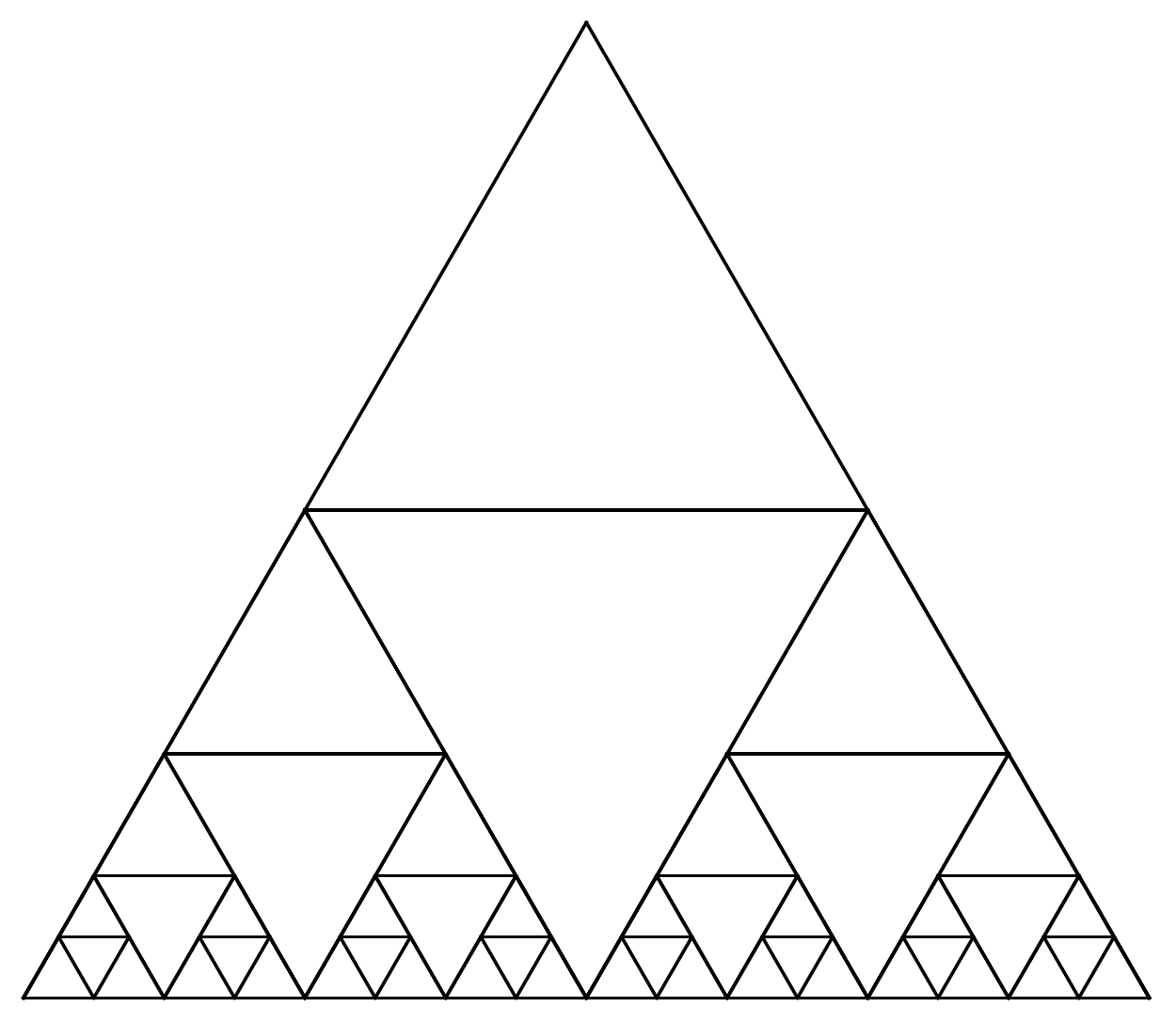}
	\begin{picture}(0,0)
	\put(-113,147){$x_{(0,1)}=q_0$}
	\put(-150,76){$x_{(1,1)}$}
	\put(-49,76){$x_{(1,2)}$}
	\put(-171,40){$x_{(2,1)}$}
	\put(-110,44){$x_{(2,2)}$}
	\put(-88,44){$x_{(2,3)}$}
	\put(-27,40){$x_{(2,4)}$}
	\end{picture}
	\caption{The points $x_{(n,k)}=F_{w(n,k)}q_0$}.\label{fig1}
\end{figure}
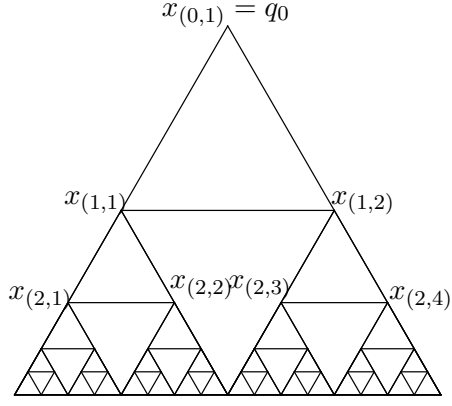

\begin{notation}	
	(a). Let $U_n=V_n\cap I$, $\tilde{U}_n=\tilde{V}_n\cap I$ and $U_*=V_*\cap I$. Clearly, $U_*$ is the set of dyadic rationals on $I$.
	 
	(b). For each pair $(n,k)$ where $n\geq 0$ and $1\leq k\leq 2^n$, let $w(n,k)$ be the unique word in $\{1,2\}^n$ such that 
	\[[\frac{k-1}{2^n},\frac{k}{2^n}]=F_{w(n,k)}([0,1]).\]
	For example, $F_{w(n,1)}=11\cdots 1$, and $F_{w(n,2^n)}=22\cdots2$.
	
	(c). Let $x_{(n,k)}=F_{w(n,k)}q_0$, and define $NU_n=\{x_{(n,k)}:1\leq k\leq 2^n\}$. See Fiugre \ref{fig1} for an illustration.
\end{notation}

With the above definitions and notations, we introduce the following space along with the extension map.

\begin{definition}\label{def36}
	(a). Let $f\in C(I)$. For $n\geq 1$ and $1\leq k\leq 2^n$, define $c_{x_{(n,k)}}=5\tilde{D}f(n,k)$. Define the extension map $\tilde{E}$ as follows,  
	$$\tilde{E}f=h+\sum_{n=1}^\infty\sum_{x\in NU_n}c_x\varphi_x,$$ where $h$ is the unique harmonic function on $\SG$ such that $h(0)=f(0),h(1)=f(1)$ and $h(\frac{1}{2})=f(\frac{1}{2})$.
	
	(b). For $\sigma>b_1$, define the space of functions on $I$
	$$\tilde{\mathcal{T}}_{\sigma}=\{f\in C(I): \sum_{n=1}^\infty\sum_{k=1}^{2^n} 5^{n\sigma}3^{-n}|\tilde{D}f(n,k)|^2<\infty\}$$
	with norm $\|f\|_{\tilde{\mathcal{T}}_{\sigma}}=\big(\|f\|^2_{L^2(I)}+\sum_{n=1}^\infty\sum_{k=1}^{2^n} 5^{n\sigma}3^{-n}|\tilde{D}f(n,k)|^2\big)^{1/2}$.
\end{definition}

Immediately from the definition, we have the following proposition. 

\begin{proposition}\label{prop37} 
   Let $b_1<\sigma<2-\frac{\log3}{\log5}$. We have $L^2_\sigma(\SG)|_I=\tilde{\mathcal{T}}_\sigma$, and $\tilde{E}$ is a continuous map from $\tilde{\mathcal{T}}_\sigma$ to $L^2_\sigma(\SG)$ such that $(\tilde{E}f)|_I=f$. 
\end{proposition}
\textit{Proof.} Let $u\in L^2_\sigma(\SG)$ with the unique expansion $u=h+\sum_{m=1}^\infty\sum_{x\in \tilde{V}_m} c_x\varphi_x$ as shown in Theorem \ref{th22}. As in the proof of Theorem \ref{thm34}, denote $\psi_0=h$ and $\psi_m=\sum_{x\in \tilde{V}_m} c_x\varphi_x$. Clearly, for $m>n+1$, we have $\tilde{D}(\psi_m|_I)(n,k)=0$ for any $1\leq k\leq 2^n$, as $\psi_m|_{U_{n+1}}=0$. In addition, for $m<n$, $\tilde{D}(\psi_m|_I)(n,k)=0$ by Lemma \ref{lemma21}. As a consequence, we have 
\[\tilde{D}(u|_I)(n,k)=\tilde{D}(\psi_n|_I)(n,k)+\tilde{D}(\psi_{n+1}|_I)(n,k).\]
Thus, 
\[\begin{aligned}
\sum_{k=1}^{2^n}|\tilde{D}(u|_I)(n,k)|^2&\leq 2\sum_{k=1}^{2^n}|\tilde{D}(\psi_n|_I)(n,k)|^2+2\sum_{k=1}^{2^n}|\tilde{D}(\psi_{n+1}|_I)(n,k)|^2\\
&\lesssim \sum_{x\in \tilde{V}_n}|c_x|^2+\sum_{x\in \tilde{V}_{n+1}}|c_x|^2.
\end{aligned}\]
Summing over the above estimate, we get $u|_I\in \tilde{\mathcal{T}}_\sigma$. Obviously $\|u|_I\|_{L^2(I)}\leq \|u|_I\|_{L^\infty(I)}\lesssim \|u\|_{L^2_\sigma(\SG)}$, so the restriction map is continuous.

Next, we show $\tilde{E}$ is the desired extension map. It is not hard to see that 
\[\tilde{D}\varphi_{x_{(n',k')}}(n,k)=\frac{1}{5}\delta_{n'n}\delta_{k'k},\]
where $\delta_{ij}$ denotes the Kronecker delta. As a consequence, 
\[\tilde{D}\big((\tilde{E}f)|_I\big)(n,k)=\frac{1}{5}c_{x_{(n,k)}}=\tilde{D}f(n,k).\]
In addition, $\tilde{E}f(0)=h(0)=f(0), \tilde{E}f(\frac{1}{2})=h(\frac{1}{2})=f(\frac{1}{2})$ and $\tilde{E}f(1)=h(1)=f(1)$. Combining the above observations, we conclude $(\tilde{E}f)|_I=f$. It is easy to check the continuity of $\tilde{E}$ with Theorem \ref{th22}. \hfill$\square$\vspace{0.2cm}

The following lemma shows the relationship between two spaces $\tilde{\mathcal{T}}_\sigma$ and $B^{2,2}_\alpha(I)$. 

\begin{lemma}\label{lemma38}
	Let $\alpha=\alpha(\sigma)$. For $b_1<\sigma<\frac{\log6}{\log5}$, we have $B^{2,2}_\alpha(I)\subset \tilde{\mathcal{T}}_\sigma$; for $b_1<\sigma<b_2$, we have $\tilde{\mathcal{T}}_\sigma=B^{2,2}_\alpha(I)$. 
\end{lemma}
\textit{Proof.} It is clear that $B^{2,2}_\alpha(I)\subset \tilde{\mathcal{T}}_\sigma$ for $b_1<\sigma<\frac{\log6}{\log5}$, as 
\[\tilde{D}f(n,k)=\begin{cases}
\frac{1}{5}A_{n+1}(f)_{2k-1}-\frac{4}{5}A_{n+1}(f)_{2k}+\frac{3}{25}A_n(f)_k+\frac{3}{25}A_n(f)_{k+1},\text{ if $k$ is odd},\\
\frac{4}{5}A_{n+1}(f)_{2k-1}-\frac{1}{5}A_{n+1}(f)_{2k}-\frac{3}{25}A_n(f)_{k-1}-\frac{3}{25}A_n(f)_k,\text{ if $k$ is even}.
\end{cases}\]
On the other hand, by Theorem \ref{thm34} and Proposition \ref{prop37}, we have $\tilde{\mathcal{T}}_\sigma=L^2_\sigma(SG)|_{I}\subset B^{2,2}_\alpha(I)$ for $b_1<\sigma<b_2$.\hfill$\square$\vspace{0.1cm}

\textbf{Remark. } One can check that the linear function $f(t)=t$ on $I$ is not in $\tilde{\mathcal{T}}_\sigma$ for $\sigma\geq \frac{\log6}{\log5}$. So the bound for the range of $\sigma$ in Lemma 3.8 is sharp. \vspace{0.2cm}

Combining Proposition \ref{prop37} and Lemma \ref{lemma38}, we get the extension theorem as follows.

\begin{theorem}\label{thm39}
Let $b_1<\sigma<b_2$ and $\alpha=\alpha(\sigma)$. The extension map $\tilde{E}$ is a continuous map from $B^{2,2}_\alpha(I)$ to $L^2_\sigma(\SG)$ such that $(\tilde{E}f)|_I=f$.
\end{theorem}

\section{A trace theorem for higher order}
In Section 3, we developed an extension of A. Jonsson's theorem. However, for Sobolev spaces of higher orders, the Besov spaces are no longer the trace spaces. In this section, we work on a discrete characterization of $L^2_\sigma(\SG)|_I$ for $b_1<\sigma\leq 2$. This includes  $L^2_2(\SG)=dom_{L^2}\Delta(SG)$ as a special case. We still study the restriction theorem and the extension theorem seperately, and prove theorem \ref{thm14} at the end.

\subsection{A restriction theorem} 
In this subsection, we will introduce the trace space (see Definition \ref{def45}) and prove a restriction theorem. 

We would like to study the space $\tilde{\mathcal{T}}_\sigma$ first, and try to modify it.

\begin{figure}[h]
	\includegraphics[width=5.2cm]{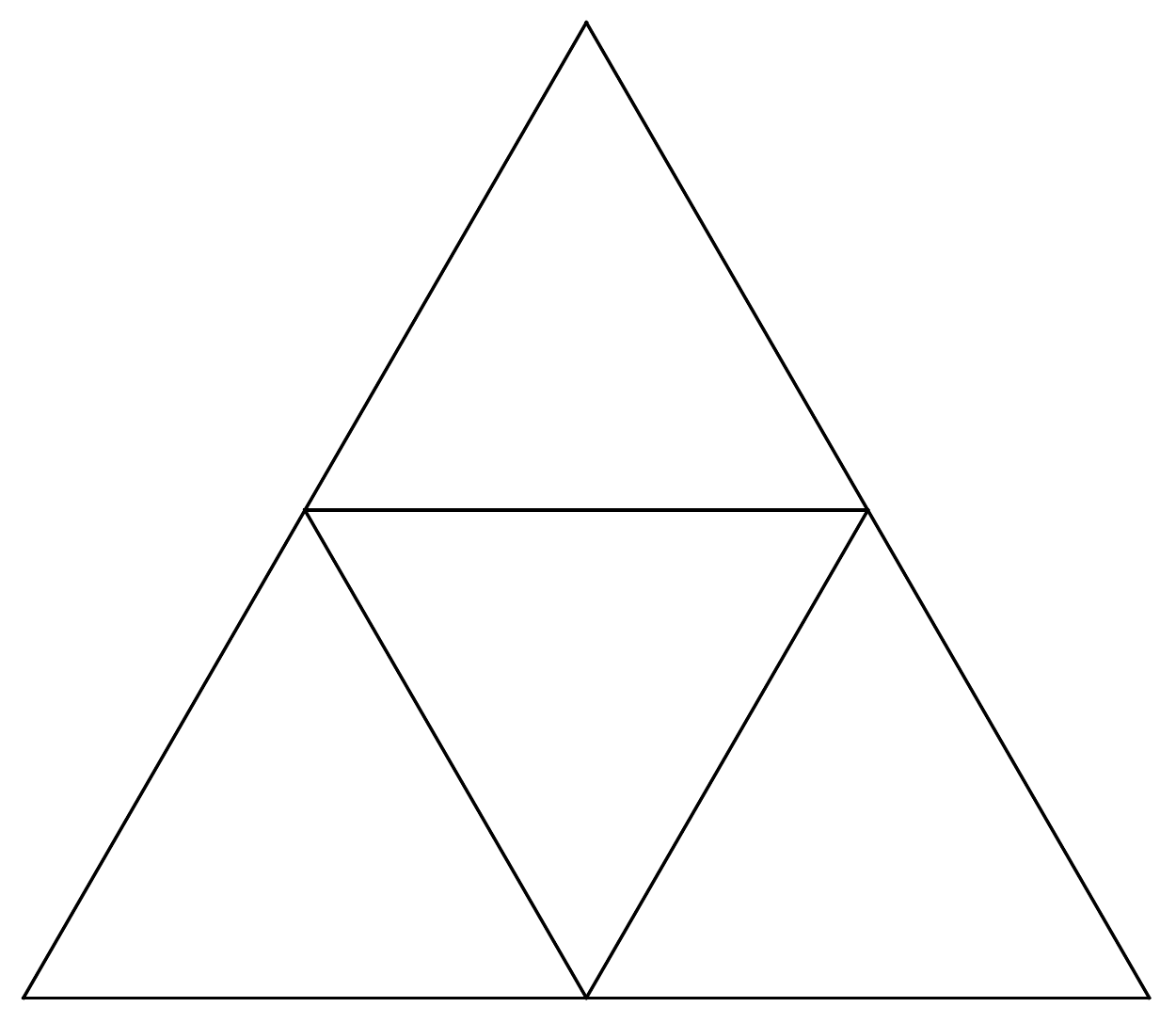}\quad
	\includegraphics[width=5.2cm]{figure1.pdf}	
	\begin{picture}(0,0)
	\put(-280,28){$Z_{1,1}$}
	\put(-210,28){$Z_{1,2}$}
	
	\put(-122,41){$\tilde{Z}_{1,1}$}
	\put(-51,41){$\tilde{Z}_{1,2}$}
	\end{picture}
	\caption{An illustration for $Z_{n,k}$ and $\tilde{Z}_{n,k}$.}\label{fig2}
\end{figure}

\begin{notation} Recall that we define $w(n,k)\in \{1,2\}^n$ such that $F_{w(n,k)}([0,1])=[\frac{k-1}{2^n},\frac{k}{2^n}]$ for $n\geq 0$ and $1\leq k\leq 2^n$.
	
(a). Write $Z_{(n,k)}=F_{w(n,k)}SG$ and $\tilde{Z}_{(n,k)}=F_{w(n,k)}F_0SG$. See Figure \ref{fig2} for an illustration of $Z_{(n,k)}$ and $\tilde{Z}_{(n,k)}$.

(b). Say $(n',k')\geq (n,k)$ if and only if $F_{w(n',k')}(I)\subset F_{w(n,k)}(I)$. It is easy to see that $$Z_{(n,k)}=\bigcup_{(n',k')\geq (n,k)}\tilde{Z}_{(n',k')}.$$ 

(c). Define non-abelian `+' on the pairs with the following equation
\[w\big((n,k)+(n',k')\big)=w(n,k)w(n',k')\in \{1,2\}^{n+n'}.\]
Clearly, $(n,k)+(n',k')\geq(n,k)$. 

As an example of (c), readers can check that 
\[(n,k)+(1,1)=(n+1,2k-1) \text{ and } (n,k)+(1,2)=(n+1,2k).\]
\end{notation}

The idea of the following lemma and Lemma \ref{lemma47} can be found in \cite{cq3}, where pointwise approximations of Laplacians were discussed. 
\begin{lemma}\label{lemma41}
	There exist $J,J'\in C(\SG)$ such that for each $u\in dom_{L^2}\Delta(SG)$ and $n\geq 0,1\leq k\leq 2^n$, we have
	\[\begin{cases}
	\tilde{D}(u|_I)(n+1,2k-1)=(\frac{3}{5})^n\int_{Z_{(n,k)}}\Delta f(x)J\circ F^{-1}_{w(n,k)}(x)d\mu(x),\\
	\tilde{D}(u|_I)(n+1,2k)=(\frac{3}{5})^n\int_{Z_{(n,k)}}  \Delta f(x)J'\circ F^{-1}_{w(n,k)}(x)d\mu(x).
	\end{cases}\]
\end{lemma}

\textit{Proof.}  First, by the Riesz representation theorem on Hilbert spaces, we can find $J\in dom_0\mathcal{E}$ such that $\tilde{D}v(1,1)=-\mathcal{E}(v,J)$ for each $v\in dom_0\mathcal{E}$. 

Define $u_0=u-h$, where $h$ is the unique harmonic function that $h|_{V_0}=u|_{V_0}$. Using the weak formula of the Laplacian and the fact that $\tilde{D}(h|_I)(1,1)=0$, we get the following desired formula
\[\begin{aligned}
\tilde{D}(u|_I)(1,1)=\tilde{D}(u_0|_I)(1,1)=-\mathcal{E}(u_0,J)=\int_{\SG} \Delta u_0\cdot Jd\mu=\int_{\SG} \Delta u\cdot Jd\mu.
\end{aligned}\] 
A same idea works for $\tilde{D}u(1,2)$. The lemma then follows by scaling.\hfill$\square$

\begin{lemma}\label{lemma42}
	The restriction map $R:L^2_\sigma(K)\to \tilde{\mathcal{T}}_{\sigma}$ is a continuous linear map for $b_1<\sigma\leq 2$. 
\end{lemma}
\textit{Proof.} It suffices to prove the argument for $1\leq \sigma\leq 2$. First, for $\sigma=1$, it is an immediate consequence of Theorem \ref{thm34} and Lemma \ref{lemma38}.

Next, we show the lemma for $\sigma=2$. By using Lemma \ref{lemma41}, we get 
\[|\tilde{D}(u|_I)(n+1,2k-1)|\lesssim (\frac{3}{5})^n\int_{Z_{(n,k)}}|\Delta u|d\mu\\
\lesssim (\frac{3}{5})^n\sum_{(n',k')\geq (n,k)}3^{-n'/2}\|\Delta u\|_{L^2(\tilde{Z}_{(n',k')})}.
\] 
The same estimate holds for $\tilde{D}(u|_I)(n+1,2k)$.  Using the above estimate and the Minkowski inequality, we have
\[\begin{aligned}
&\big(\sum_{n=1}^\infty\sum_{k=1}^{2^n} 5^{2n}\cdot 3^{-n}|\tilde{D}(u|_I)(n,k)|^2\big)^{1/2}\\
\lesssim& \big\|3^{n/2}\sum_{(n',k')\geq (n,k)}3^{-n'/2}\|\Delta u\|_{L^2(\tilde{Z}_{(n',k')})}\big\|_{l^2(n,k)}
\\=&\big\|\sum_{(n',k')\geq (0,1)}3^{-n'/2}\|\Delta u\|_{L^2(\tilde{Z}_{(n,k)+(n',k')})}\big\|_{l^2(n,k)}\\
\leq& \sum_{(n',k')\geq (0,1)}3^{-n'/2}\|\Delta u\|_{L^2(\bigcup_{(n,k)\geq (0,1)}\tilde{Z}_{(n,k)+(n',k')})}\\
\leq &\sum_{n'=0}^\infty (\frac{2}{3})^{n'/2} \|\Delta u\|_{L^2(\bigcup_{k=1}^{2^{n'}}Z_{(n',k)})}\lesssim \|\Delta u\|_{L^2(\SG)},
\end{aligned}\]
where we use the notation $\|a_{(n,k)}\|_{l^2(n,k)}=\big(\sum_{(n,k)\geq (0,1)}|a_{(n,k)}|^2\big)^{1/2}$ for convenience. This proves the argument for $\sigma=2$. 

For general $1\leq \sigma\leq 2$, we can use complex interpolation to deduce the lemma.\hfill$\square$\vspace{0.2cm}

However, $\tilde{\mathcal{T}}_{2}$ is actually a larger space than the trace space of $L^2_2(SG)$. For example, we will see in Corollary \ref{coro33} that $\varphi_x|_{I}\in \tilde{\mathcal{T}}_2\setminus \big(L^2_2(SG)|_I\big)$, where $\varphi_x$ is a tent function.

Recall that on $\SG$, for a function $u\in C(\SG)$, we define the normal derivative at a boundary point to be $$\partial_n u(q_i)=\lim_{n\to\infty}(\frac53)^n \big(2u(q_i)-u(F_i^nq_j)-u(F_i^nq_k)\big),\text{ with }\{i,j,k\}=\{0,1,2\}.$$
The definition can be localized to any vertex in $V_*$ by scaling, and we use $\uparrow,\rightarrow,\leftarrow$ to show the direction, i.e.
\[\begin{aligned}
\partial^\uparrow_n u(F_wq_0)=(\frac{5}{3})&^{|w|}\partial_n (u\circ F_w)(q_0),\\ \partial^\leftarrow_n u(F_wq_1)=(\frac{5}{3})^{|w|}\partial_n(u\circ F_w)(q_1),&\quad \partial^\rightarrow_n u(F_wq_2)=(\frac{5}{3})^{|w|}\partial_n(u\circ F_w)(q_2).
\end{aligned}\]  

\begin{lemma}\label{lemma43} 
	Let $u\in dom_{L^2}\Delta(\SG)$. We have 
	 $$\partial_n u(q_1)=\lim_{n\to\infty}(\frac{5}{3})^n\big(4u(q_1)-5u(F_1^{n+1}q_2)+u(F^n_1q_2)\big).$$
	For general cases, for $n\geq 0$ and $1\leq k\leq 2^n$,
	$$\partial^\rightarrow_nu(\frac{k}{2^n})=\lim_{m\to\infty}(\frac{5}{3})^m\big(4u(\frac{k}{2^n})-5u(\frac{k}{2^n}-\frac{1}{2^{m+1}})+u(\frac{k}{2^n}-\frac{1}{2^{m}})\big);$$
	for $n\geq 0$ and $0\leq k\leq 2^n-1$
	$$\partial^\leftarrow_nu(\frac{k}{2^n})=\lim_{m\to\infty}(\frac{5}{3})^m\big(4u(\frac{k}{2^n})-5u(\frac{k}{2^n}+\frac{1}{2^{m+1}})+u(\frac{k}{2^n}+\frac{1}{2^{m}})\big).$$
\end{lemma}

\textit{Proof.} We only need to show the special case for $q_1=0$, since general cases can be proven by using scaling and symmetry.

First, the equation holds for harmonic functions without taking the limit, since 
\[u(F_1^{n+1}q_2)=\frac{1}{5}u(F^n_1q_0)+\frac{2}{5}u(q_1)+\frac{2}{5}u(F^n_1q_2).\]
For general $u\in  dom_{L^2}\Delta(\SG)$, we only need to notice that 
\[u(F_1^{n+1}q_2)=\frac{1}{5}u(F^n_1q_0)+\frac{2}{5}u(q_1)+\frac{2}{5}u(F^n_1q_2)+\frac{1}{5^n}\int_{SG}G(F_1q_2,y)\Delta u(F_1^ny)d\mu(y),\]
where $G$ is the Green's function on $\SG$. \hfill$\square$

\begin{corollary}
	Let $x=\frac{1}{2}$, then $R\varphi_x\notin RL^2_2(SG)$.
\end{corollary}
\textit{Proof}. Assume there exists $u\in dom_{L^2}\Delta(\SG)$ such that $u|_{I}=\varphi_x|_{I}$. Then by Lemma \ref{lemma43}, we have 
\[\partial^\rightarrow_n u(\frac{1}{2})=\partial^\rightarrow_n \varphi_x(\frac{1}{2}),\quad \partial^\leftarrow_n u(\frac{1}{2})=\partial^\leftarrow_n \varphi_x(\frac{1}{2}).\]
Thus $u$ does not satisfies the matching condition at $x$, i.e. $\partial^\rightarrow_n u(\frac{1}{2})+\partial^\leftarrow_n f(\frac{1}{2})\neq 0$, which contradicts the fact that $u\in dom_{L^2}\Delta(\SG)$.\hfill$\square$

Inspired by the above observation, we need to include the information of matching condition into the desired trace space. 

\begin{definition}\label{def45}	Let $f\in C(I)$, and let $n\geq 2,1\leq k\leq 2^n-1$.
	
	(a). Define $Df(n,k),$ as follows. For $k$ odd, define 
	\[Df(n,k)=\tilde{D}f(n-1,\frac{k+1}{2});\]
	for $k$ even, define 
	\[Df(n,k)=f(\frac{k}{2^n})-\frac{5}{8}\big(f(\frac{k-1}{2^n})+f(\frac{k+1}{2^n})\big)+\frac{1}{8}\big(f(\frac{k-2}{2^n})+f(\frac{k+2}{2^n})\big).\]
	
	(b). Define
	\[
	\mathcal{T}_{\sigma}=\big\{f\in C(I):\sum_{n=2}^\infty\sum_{k=1}^{2^n-1} 5^{\sigma n}3^{-n}|Df(n,k)|^2<\infty\big\},
	\]
	with norm $\|f\|_{\mathcal{T}_\sigma}=\big(\|f\|^2_{L^2(I)}+\sum_{n=2}^\infty\sum_{k=1}^{2^n-1} 5^{\sigma n}3^{-n}|Df(n,k)|^2\big)^{1/2}$.
\end{definition}

\textbf{Remark.} We can also characterize $\mathcal{T}_\sigma$ with 
\[\mathcal{T}_\sigma=\{f\in\tilde{\mathcal{T}}_\sigma:\sum_{n=2}^\infty\sum_{k=1}^{2^{n-1}-1}5^{\sigma n}3^{-n}|Df(n,2k)|^2<\infty\},\]
which means we additionally require the matching condition on $\mathcal{\tilde{T}}_\sigma$. In addition, for small $\sigma$, the two spaces coincide as stated by the following lemma. 

\begin{lemma}\label{lemma46}
	For $b_1<\sigma<b_2$, we have $\mathcal{T}_\sigma=\tilde{\mathcal{T}}_\sigma$.
\end{lemma}
\textit{Proof.} By the above remark and using Lemma \ref{lemma38}, we can easily check $\mathcal{T}_\sigma\subset \tilde{\mathcal{T}}_\sigma=B^{2,2}_\alpha(I)\subset \mathcal{T}_\sigma$, where $\alpha=\alpha(\sigma)$.\hfill$\square$\vspace{0.2cm}

Lemma \ref{lemma46} can be polished, see Corollary \ref{coro411}. Parellel to Lemma \ref{lemma41}, we have the following lemma \ref{lemma47}.

\begin{lemma}\label{lemma47}
	There exists $J_1,J_2\in C(\SG)$ such that for each $u\in dom_{L^2}\Delta(\SG)$ and $n\geq 1,1\leq k\leq 2^n-1$, the following equality holds 
	\[\begin{aligned}
    D(u|_I)(n+1,2k)=&(\frac{3}{5})^n\int_{F_{w(n,k)}SG}J_1(F^{-1}_{w(n,k)}x)\Delta u(x)d\mu(x)\\&+(\frac{3}{5})^n\int_{F_{w(n,k+1)}SG}J_2(F^{-1}_{w(n,k+1)}x)\Delta u(x)d\mu(x). 
    \end{aligned}\]
\end{lemma} 
\textit{Proof.}  The proof is very similar to that of Lemma \ref{lemma41}. For any function $h$ that is harmonic on $F_1\SG\cup F_2\SG$, it is direct to check that $D(h|_I)(2,2)=0$. Let $\tilde{u}=u-h$ on $F_1\SG\cup F_2\SG$, where $h$ is harmonic on $F_1\SG\cup F_2\SG$ with boundary values $h(q_1)=u(q_1),h(F_0q_1)=u(F_0q_1),h(F_0q_2)=u(F_0q_2)$ and $h(q_2)=u(q_2)$. We can find $J_0$ on $F_1\SG\cup F_2\SG$ such that $\mathcal{E}_{F_1\SG\cup F_2\SG}(\tilde{u},J)=-D(\tilde{u}|_I)(2,2)$. By a same argument as in the proof of Lemma \ref{lemma41}, we have
\[D(u|_I)(2,2)=\int_{F_1\SG\cup F_2\SG} J_0\Delta ud\mu.\]
Take $J_1=\frac{5}{3}J_0\circ F_1$ and $J_2=\frac{5}{3}J_0\circ F_2$, then we get the desired equation for $D(u|_I)(2,2)$. For general cases, we only need to use scaling.\hfill$\square$\vspace{0.2cm}

Following a same proof of Lemma \ref{lemma42}, we finally get the following restriction theorem.

\begin{theorem}\label{thm48}
	The restriction map $R:L^2_\sigma(K)\to \mathcal{T}_\sigma$ is a continuous linear map for $b_1<\sigma\leq 2$. 
\end{theorem}
\textit{Proof. } For $b_1<\sigma<b_2$, the result is an easy consequence of Lemma \ref{lemma46} and Proposition \ref{prop37}.  For $\sigma=2$, the result follows from a similar estimate as the proof of Lemma \ref{lemma42}, using Lemma \ref{lemma47} and Lemma \ref{lemma42}. The theorem then follows by using complex interpolation.\hfill$\square$

\subsection{An extension theorem}
Next, we construct an extension map by modifying $\tilde{E}$. 

\begin{definition}\label{extension}
	Choose $v_0\in dom_{L^2}\Delta(\SG)$ such that $v_0|_{V_0}=0$ and $\partial_n v_0(q_j)=\delta_{0,j}$. Recall that in Definition \ref{def36} $$\tilde{E}f=h+\sum_{n=1}^\infty\sum_{x\in NU_n} c_x\varphi_x,$$ where $c_{x_{(n,k)}}=5\tilde{D}f(n,k)=5Df(n+1,2k-1)$ and $h\in \mathcal{H}_0$ with $h|_{U_1}=f|_{U_1}$. Define the extension map $E$ as follows,
	$$Ef=\tilde{E}f-\sum_{n=1}^\infty\sum_{k=1}^{2^n}\frac{12}{5}c_{x_{(n,k)}}v_0\circ F^{-1}_0\circ F^{-1}_{w(n,k)}.$$
\end{definition}

\begin{theorem}\label{thm410}
For $b_1<\sigma\leq2$, the extension map $E:\mathcal{T}_\sigma\to L^2_\sigma(SG)$ is a continuous linear map such that $(Ef)|_{I}=f$.
\end{theorem}
\textit{Proof.} First, we show $E:\mathcal{T}_\sigma\to L^2_\sigma(\SG)$ is bounded for $1\leq\sigma\leq2$. Let 
$$\tilde{E}_mf=h+\sum_{n=1}^m\sum_{x\in NU_n} c_x\varphi_x.$$
Choose $v_2\in dom\Delta$ such that $v_2|_{V_0}=0$ and $\partial_n v_2(q_j)=\delta_{2,j}$, and define 
\[\begin{aligned}
E_mf(x)=&\tilde{E}_mf(x)-\sum_{n=1}^m\sum_{k=1}^{2^n}\frac{12}{5}c_{x_{n,k}}v_0\circ F^{-1}_0\circ F^{-1}_{w(n,k)}\\
&-\sum_{k=1}^{2^m-1}\big(\frac{24}{5}Df(m+1,2k)\big)v_2\circ F^{-1}_{w(m+1,2k)}.
\end{aligned}\] 
$E_mf$ satisfies the matching conditions at all vertices, which implies that $E_mf\in L^2_2(\SG)\subset L^2_\sigma(\SG)$.

Notice that $v_0\circ F^{-1}_0\circ F^{-1}_{w(n,k)}$ supports on $\tilde{Z}_{(n,k)}$, and $v_2\circ F^{-1}_{w(m+1,2k)}$ supports on $F_{w(m+1,2k)}SG$, which are disjoint sets. We can easily get the following estimates
$$\begin{cases}
\|E_mf\|_{L^2_1}\lesssim \big(f^2(0)+f^2(1)+f^2(\frac{1}{2})+\sum_{n=1}^m\sum_{k=1}^{2^n}(\frac{5}{3})^n|\tilde{D}f(n,k)|^2\\\qquad\qquad\qquad+\sum_{k=1}^{2^m-1}(\frac{5}{3})^m|Df(m+1,2k)|^2\big)^{1/2},\\
\|E_mf\|_{L^2_2}\lesssim \big(f^2(0)+f^2(1)+f^2(\frac{1}{2})+\sum_{n=1}^m\sum_{k=1}^{2^n}(\frac{25}{3})^n|\tilde{D}f(n,k)|^2\\\qquad\qquad\qquad+\sum_{k=1}^{2^m-1}(\frac{25}{3})^m|Df(m+1,2k)|^2\big)^{1/2}.
\end{cases}$$
Clearly, the above estimates holds uniformly for any $m\geq 1$. Using complex interpolation, we then get 
\[
\begin{aligned}
\|E_mf\|_{L^2_\sigma}\lesssim \big(f^2(0)+f^2(1)+f^2(\frac{1}{2})+\sum_{n=1}^m\sum_{k=1}^{2^n}5^{\sigma n}3^{-n}|\tilde{D}f(n,k)|^2\\\qquad\qquad\qquad+\sum_{k=1}^{2^m-1}5^{\sigma m}3^{-m}|Df(m+1,2k)|^2\big)^{1/2}
\end{aligned}
\]
holds uniformly for any $m$. In other words, $\|E_mf\|_{L^2_\sigma}\lesssim \|f\|_{\mathcal{T}_\sigma}$ uniformly for any $m\geq 1$.

Similarly, the following estimate holds uniformly for any $1\leq\sigma\leq2$ and $m\geq 1$,
\[\begin{aligned}
\|E_{m'}f-E_mf\|^2_{L^2_\sigma(SG)}&\lesssim \sum_{n=m+1}^\infty\sum_{k=1}^{2^n}5^{2n}3^{-n}|\tilde{D}f(n,k)|^2+\sum_{k=1}^{2^{m}-1}5^{2m}3^{-m}\big(Df(m+1,2k)\big)^2\\
&+\sum_{k=1}^{2^{m'}-1}5^{2m'}3^{-m'}\big(Df(m'+1,2k)\big)^2.
\end{aligned}\]
As a result, $E_mf$ converges in $L^2_\sigma(\SG)$. Noticing that $E_mf$ converges pointwise to $Ef$, we conclude that $E_mf$ converges to $Ef$ in $L^2_\sigma(\SG)$ sense. Thus $Ef\in L^2_\sigma(\SG)$ and $\|Ef\|_{L^2_\sigma(\SG)}\lesssim \|f\|_{\mathcal{T}_\sigma}$.

Next, for $b_1<\sigma<1$, we have the scaling property that $\|v_0\circ F_w^{-1}\|_{L^2_\sigma(\SG)}\asymp (5^\sigma3^{-1/2})^{|w|}\|v_0\|_{L^2_\sigma(\SG)}$, as a consequence of Theorem \ref{th22}. In addition, using Theorem \ref{th22}, we can check that
\[\|\sum_{n=1}^\infty\sum_{k=1}^{2^n}\frac{12}{5}c_{x_{(n,k)}}v_0\circ F^{-1}_0\circ F^{-1}_{w(n,k)}\|^2_{L^2_\sigma(\SG)}\asymp \sum_{n=1}^\infty\sum_{k=1}^{2^n}\|\frac{12}{5}c_{x_{(n,k)}}v_0\circ F^{-1}_0\circ F^{-1}_{w(n,k)}\|^2_{L^2_\sigma(\SG)},\]
as $v_0\circ F^{-1}_0\circ F^{-1}_{w(n,k)}$ have disjoint supports. Combining the above two facts, it is direct to see that $E:\mathcal{T}_\sigma\to L^2_\sigma(\SG)$ is continuous for small $\sigma$.

Lastly, $(Ef)|_{I}=(\tilde{E}f)|_I=f$, since $v_0\circ F_0^{-1}\circ F_{w(n,k)}^{-1}$ is supported away from $I$. \hfill$\square$

Combining Theorem \ref{thm48} and \ref{thm410}, we finally get Theorem \ref{thm14}. Also, the following corollary shows the relationship of the different traces spaces discussed in this paper.

\begin{corollary}\label{coro411} Let $\alpha=\alpha(\sigma)$. Then we have the following relationships.
	
	(a) For $b_1<\sigma<b_2$, we have $B^{2,2}_\alpha(I)=\mathcal{T}_\sigma$; for $b_1<\sigma<\frac{\log6}{\log5}$, we have $B^{2,2}_\alpha(I)=\mathcal{T}_\sigma$; for $b_1<\sigma<2-\frac{\log3}{\log5}$, we have $\tilde{\mathcal{T}}_\sigma=\mathcal{T}_\sigma$.
	
	(b) For any $b_1<\sigma_1<\sigma_2\leq 2$, we have $[\mathcal{T}_{\sigma_1},\mathcal{T}_{\sigma_2}]_\theta=\mathcal{T}_{(1-\theta)\sigma_1+\theta\sigma_2}$.
\end{corollary}
\textit{Proof.} (a) is immediately from Lemma \ref{lemma38}, Proposition \ref{prop37} and Theorem \ref{thm14}.

(b). Using complex interpolation, the restriction map maps from $L^2_{(1-\theta)\sigma_1+\theta\sigma_2}(\SG)$
to $[\mathcal{T}_{\sigma_1},\mathcal{T}_{\sigma_2}]_\theta$, and the extension map maps from $[\mathcal{T}_{\sigma_1},\mathcal{T}_{\sigma_2}]_\theta$ to $L^2_{(1-\theta)\sigma_1+\theta\sigma_2}(\SG)$. Thus $\mathcal{T}_{(1-\theta)\sigma_1+\theta\sigma_2}=L^2_{(1-\theta)\sigma_1+\theta\sigma_2}(\SG)|_I=[\mathcal{T}_{\sigma_1},\mathcal{T}_{\sigma_2}]_\theta$. \hfill$\square$

\section{Related observations and further questions}
In this last section, we provide some related results and question that worth further study. 

Another space that we are interested in is  
$$dom_{L^\infty}\Delta(\SG)=\{u\in dom\mathcal{E}:\Delta u\in L^\infty(\SG)\}.$$
With a same method as in the last section, we can derive the following result.
 
\begin{theorem}\label{thm51}
Define
$$\mathcal{T}^\infty_2=\{f\in C(\SG):\sup_{n\geq 2}\sup_{1\leq k\leq 2^n-1}5^n|Df(n,k)|<\infty\}.$$
Then, $dom_{L^\infty}\Delta(\SG)|_{I}=\mathcal{T}^\infty_2$.
\end{theorem}

Consider the symmetric derivative of the functions. Let $f\in C(I)$ and fix $x\in U_m\setminus U_0$, we define $\delta_nf(x)=f(x-\frac{1}{2^n})-2f(x)+f(x+\frac{1}{2^n}).$
The symmetric derivative at $x$ is defined to be the renormalized limit of $\delta_n f(x)$, 
\[f''_s(x)=\lim_{n\to\infty} 4^n\delta_nf(x).\]

\begin{proposition}\label{prop51}
	Let $u\in dom_{L^\infty}\Delta(\SG)$ and $f=u|_I$. Then for all $x\in U_*\setminus \{0,1\}$, we have $g''_s(x)=0$.
\end{proposition}

\textit{Proof.} Let $x=\frac{k}{2^n}$ for some $n\geq 1$ and $1\leq k\leq 2^n-1$. By direct computation and using Theorem \ref{thm51}, we have 
\[|\delta_{m+1}f(x)-\frac{1}{5}\delta_mf(x)|=|\frac{8}{5}Df(m+1,2^{m+1-n}k)|\leq c(\frac{1}{5})^{m-n},\]
where $c$ is a constant depends on $\|\Delta u\|_\infty$. Summing over the above estimate, we get
\[|\delta_{m}f-(\frac{1}{5})^{m-n+1}\delta_{n+1}f|\leq c(m-n-1)(\frac15)^{m-n}.\]
As an immediate consequence, we get $f''_s(x)=\lim_{m\to\infty} 4^m\delta_{m}f(x)=0$.\hfill$\square$\vspace{0.2cm}

On the other hand, $\delta_nf(x)$ should not converge uniformly to $0$ in general cases. Since otherwise, it would imply that $f\in dom\Delta (I)$ with $\Delta f=0$, which means $f$ is a linear function. In fact, the following result shows that $f'(x)=\lim\limits_{n\to\infty}2^n(f(x)-f(x+\frac{1}{2^n}))$ diverges when $\partial^\rightarrow_nu(x)\neq 0$.

\begin{proposition}
	Let $u\in dom_{L^\infty}\Delta(\SG)$ and let $f=u|_I$. Let $x\in U_*$ and suppose $\partial_n u(x)\neq 0$. Then $|f'(x)|=\infty$.
\end{proposition}

\textit{Proof.} Without loss of generality, assume $x\neq 1$, and we denote $\delta^+_mf(x)=f(x+\frac{1}{2^m})-f(x)$. Then by Lemma \ref{lemma43}, we have
\[\begin{aligned}
\lim_{m\to\infty} (\frac{5}{3})^m \big(\delta^+_mf(x)-5\delta^+_{m+1}f(x)\big)=\partial_n^\rightarrow u(x).
\end{aligned}\]
Clearly, 
\[\begin{aligned}
\delta_m^+f(x)&=\sum_{k=n}^{m-1} \big((\frac{1}{5})^{m-k-1}\delta^+_{k+1}f(x)-(\frac{1}{5})^{m-k}\delta^+_kf(x)\big)+(\frac{1}{5})^{m-n}\delta^+_n f(x)\\
&\asymp -\sum_{k=n}^{m-1} (\frac{3}{5})^{k}(\frac{1}{5})^{m-k}\big(\partial_n^\rightarrow u(x)+o(1)\big)+(\frac{1}{5})^{m-n}\delta^+_n f(x).
\end{aligned}\]
As a result, we have $\lim\limits_{m\to\infty}(\frac{5}{3})^m\delta^+_m f(x)=-\frac{1}{2}\partial^\rightarrow_n u(x)$. Similarly, for $x\neq 0$, let $\delta^-_mf(x)=f(x)-f(x-\frac{1}{2^m})$, we have $\lim\limits_{m\to\infty}(\frac{5}{3})^m\delta^-_m f(x)=\frac{1}{2}\partial^\leftarrow_n u(x)$. The proposition is immediate from the above observation. \hfill$\square$\vspace{0.2cm}

We are also interested in Sobolev spaces of higher orders. We believe that a similar idea would work, but more complicated differences will occur in the discrete characterization. For example, for $1<\sigma\leq 2$, we will need to study extension algorithm of biharmonic functions, and find suitable difference operators. The computations are getting messy, so we do not go further in this direction. Hopefully, readers may get new ideas dealing with this.

Readers may have noticed that $b_2$ plays the important role in that it is the highest index that $B^{2,2}_{\alpha}(I)\supset L^2_\sigma(\SG)|_{I}$ for each $\alpha<\alpha(b_2)$ as long as $\sigma\geq\alpha(b_2)$. Noticing that $b_2$ is uniquely characterized by harmonic functions, we wonder whether biharmonic functions play similarly important roles in higher order cases. It is also of interest to find the largest index $\alpha$ such that $h|_{I}$ lies in $B^{2,2}_{\alpha}(I)$ for multi-harmonic functions, and what kind of role these indexes will play. We hope to find a systematic way to deal with these.

\bibliographystyle{amsplain}

\end{document}